\documentclass[12pt]{amsart}
\usepackage[margin=1in]{geometry}
\usepackage{times}
\usepackage{amsthm,amssymb}
\usepackage{amsmath,xypic}
\usepackage[usenames,dvipsnames]{color}
\usepackage{hyperref}
\usepackage{epigraph}
\usepackage{graphicx}
\usepackage[authoryear]{natbib}
\usepackage{fancyhdr} 
\usepackage{enumitem}
\newcommand{\be}{\begin{equation}}
\newcommand{\ee}{\end{equation}}
\newcommand{\bes}{\begin{equation*}}
\newcommand{\ees}{\end{equation*}}
\newcommand{\bea}{\begin{eqnarray}}
\newcommand{\eea}{\end{eqnarray}}
\newcommand{\beas}{\begin{eqnarray}}
\newcommand{\eeas}{\end{eqnarray}}
\newcommand{\ben}{\begin{note}}
\newcommand{\een}{\end{note}}
\newcommand{\bexl}{\vskip0.1em\noindent\hrulefill\vskip1em\begin{ExerciseList}}
\newcommand{\eexl}{\end{ExerciseList}\hrulefill}

\newcommand{\bthm}{\begin{theorem}}
\newcommand{\ethm}{\end{theorem}}
\newcommand{\bpro}{\begin{prop}}
\newcommand{\epro}{\end{prop}}
\newcommand{\bcor}{\begin{corollary}}
\newcommand{\ecor}{\end{corollary}}
\newcommand{\bcon}{\begin{conjecture}}
\newcommand{\econ}{\end{conjecture}}
\newcommand{\bp}{\begin{proof}}
\newcommand{\ep}{\end{proof}}
\newcommand{\blem}{\begin{lemma}}
\newcommand{\elem}{\end{lemma}}
\newcommand{\bn}{\begin{note}}
\newcommand{\en}{\end{note}}
\newcommand{\benum}{\begin{enumerate}}
\newcommand{\eenum}{\end{enumerate}}
\newcommand{\bed}{\begin{defn}}
\newcommand{\eed}{\end{defn}}
\newcommand{\brem}{\begin{remark}}
\newcommand{\erem}{\end{remark}}

\newcommand{\btik}{\begin{tikzpicture}\begin{axis}[scale=0.5,axis y line=center, axis x line=middle]}
\newcommand{\etik}{\end{axis}\end{tikzpicture}}

\let\into=\hookrightarrow

\newcommand{\upperRomannumeral}[1]{\uppercase\expandafter{\romannumeral#1}}

\newtheorem{theorem}[equation]{Theorem}      
\newtheorem{lemma}[equation]{Lemma}          %
\newtheorem{corollary}[equation]{Corollary}  
\newtheorem{proposition}[equation]{Proposition}

\theoremstyle{definition}

\theoremstyle{definition}
\newtheorem{defn}[equation]{Definition}
\theoremstyle{remark}

\theoremstyle{definition}
\newtheorem{remark}[equation]{Remark}

\numberwithin{equation}{section}

\let\into=\hookrightarrow
\let\isom=\simeq

\let\tensor=\otimes

\newcommand{\Z}{{\mathbb Z}}

\renewcommand{\int}{\operatorname{int}}
\renewcommand{\O}{{\mathcal O}}
\renewcommand{\P}{{\mathbb P}}

\renewcommand{\bpro}{\begin{proposition}}
\renewcommand{\epro}{\end{proposition}}

\let\citep=\cite 
\renewcommand{\cite}{\citeauthor*}

\begin{document}

\title[]{On the construction of Weakly Ulrich bundles}%
\author{Kirti Joshi}%
\address{Math. department, University of Arizona, 617 N Santa Rita, Tucson
85721-0089, USA.} \email{kirti@math.arizona.edu}

\thanks{}%


\begin{abstract}
A well-known conjecture of Eisenbud, Schreyer and Weyman suggests that any projective variety carries an Ulrich (and hence also weakly Ulrich) bundle. This is known only in a handful of cases.	In this paper I provide weakly Ulrich bundles on a class of surfaces and threefolds. I construct a weakly Ulrich bundle of large rank on any smooth complete surface  in $\P^3$ over fields of characteristic $p>0$ and also for some classes of surfaces of general type in $\P^n$. I also construct intrinsic weakly Ulrich bundles on  any Frobenius split variety of dimension at most three. The bundles constructed here are in fact ACM and weakly Ulrich bundles and so I call them almost Ulrich bundles. Frobenius split varieties in dimension three include as special cases: (1) smooth hypersurfaces in $\P^4$ of degree at most four, (2) more generally, Frobenius split Fano varieties of dimension at most three, (3) Frobenius split, smooth quintics in $\P^4$ (4) more generally, Frobenius split Calabi-Yau varieties of dimension at most three (5) Frobenius split (i.e ordinary) abelian varieties of dimension at most three.   These results also imply that Chow form of these varieties  is the support of a single intrinsic determinantal equation.
\end{abstract}
\maketitle
%

\newcommand{\pn}[1]{\mathbb{P}^{#1}}
\newcommand{\pnn}{\pn{n}}
\newcommand{\pnno}{\O_{\pnn}}
\tableofcontents
\section{Introduction} 
It has been conjectured  in (\cite{eisenbud03}) that Ulrich bundles (and hence weakly Ulrich bundles) exist on any projective variety. At the moment this is known only for a handful of cases and mostly in characteristic zero  (see \cite{beauville17}).  In Theorem~\ref{th:main} I construct, under mild assumptions on the degree and (positive) characteristic,  a weakly Ulrich bundle on surfaces of general type in $\P^3$. As will be evident from my proof, a similar construction also provides weakly Ulrich bundles on smooth, projective complete intersection surfaces of general type. All of the bundles constructed here are also arithmetically Cohen-Macaulay (ACM) bundles and I call such bundles \emph{almost Ulrich bundles}.

Let me remind the reader that the importance of  weakly Ulrich (resp. Ulrich) bundles stems from the central result of (\cite{eisenbud03}) which asserts that if $X\into \P^n$ is equipped with a weakly Ulrich (resp. Ulrich) bundle for this embedding then the Chow form of $X$ (with respect to this embedding) is the support of a single determinantal equation (resp. single linear determinantal equation).

In Theorem~\ref{th:gen-type} I provide a construction of almost Ulrich bundles on a class of surfaces of general type which includes smooth quintic surfaces in $\P^3$.

In Theorem~\ref{th:calabi-yau} I  provide a construction of an intrinsic almost Ulrich bundle on any Frobenius split variety of dimension at most three. In (Corollaries~\ref{cor:calabi-yau}--\ref{cor:fano}) I record a number of special cases which are of particular interest: any Frobenius split (i.e. ordinary) abelian surface (or threefold), any Frobenius split (i.e. ordinary) K3 surface, and any ordinary or Frobenius split, Calabi-Yau threefold and finally Frobenius split Fano surface or threefold.  In particular to quadric and cubic surfaces (in $\P^3$) and quadric, cubic or quartic hypersurfaces in $\P^4$.  This also applies (rather trivially) to any embeddings $\P^2\into\P^n$ (resp. $\P^3\into\P^n$) as any projective space is Frobenius split (see \cite{joshi03} for comparison of Frobenius splitting and ordinarity).

In particular it follows from these results together with the aforementioned result  of (\cite{eisenbud03}) that the respective varieties to which the theorems apply have Chow forms which are supports of a single determinantal equation (as shown in loc. cit.). Another feature of my construction is that the bundles constructed in these theorems are intrinsic to $X$ i.e. essentially independent of the given embedding. 

Especially it follows from Theorem~\ref{th:calabi-yau}, Corollary~\ref{cor:calabi-yau}, Corollary~\ref{cor:abelian}, Corollary~\ref{cor:fano} and (\cite{eisenbud03}) that over any algebraically closed field of any positive characteristic, the Chow form of every Frobenius split Calabi-Yau, Abelian or Fano variety of dimension at most three is the support of a single intrinsic determinantal equation with respect to any projective embedding. This is the content of Theorem~\ref{th:chow-form} (for a discussion of intrinsic equations see \cite{joshi19}).

As far as I am aware, no result of comparable strength is known for Calabi-Yau, Abelian or Fano threefolds in characteristic zero.

It is a pleasure to thank N. Mohan Kumar for comments and suggestions. Thanks are due to the referee for suggestions and corrections which have improved the readability of this paper.
\section{Ulrich, weakly Ulrich and almost Ulrich bundles}
In this section the following notations are in force.
Let $X$ be a smooth, projective variety of dimension $\dim(X)=q$  with $X\into\P^n$ and equipped with $\O_X(1)$ provided by this embedding. Let $E$ be a vector bundle on $X$. Then following (\cite{eisenbud03}) one says that $E$ is a \emph{weakly Ulrich bundle} on $X$ if 
\benum 
\item $H^j(E(-m))=0$ for $1\leq j\leq q= \dim(X)$ and $m\leq j-1\leq q-1= \dim(X)-1$, and
\item  $H^j(E(-m))=0$ for $0\leq j\leq \dim(X)-1$ and $m\geq j+2$.
\eenum  
One says that $E$ is an \emph{Ulrich bundle} on $X$ if it satisfies the following conditions:
\benum
\item $H^q(E(-m))=0$ for all $m\leq q$ (here $q=\dim(X)$),
\item $H^j(E(m))=0$ for $1\leq j\leq q-1= \dim(X)-1 \text{ and all } m\in \Z$ (i.e. $E$ is arithmetically Cohen-Macaulay), and
\item $H^0(E(-m))=0$ for all $m\geq 1$.
\eenum
Reader should consult (\cite{eisenbud03}) or the excellent survey (\cite{beauville17}) for more on  Ulrich and weakly Ulrich bundles.

The bundles I construct here are not Ulrich bundles but they satisfy a stronger condition than being weakly Ulrich. I call such bundles \emph{almost Ulrich bundles}. More precisely one says that $E$ is an \emph{almost Ulrich bundle} if $E$ is weakly Ulrich and an ACM bundle.  Explicitly a bundle is \emph{almost Ulrich bundle} if and only if the following conditions are satisfies:
\benum
\item $H^q(E(-m))=0$ for all $m\leq q-1=\dim(X)-1$ (here $q=\dim(X)$),
\item $H^j(E(m))=0$ for $1\leq j\leq \dim(X)-1 \text{ and all } m\in \Z$ (i.e. $E$ is arithmetically Cohen-Macaulay), and
\item $H^0(E(-m))=0$ for all $m\geq 2$.
\eenum
In particular a bundle $E$ is Ulrich if and only if $E$ is almost Ulrich and $$H^q(E(-q))=0=H^0(E(-1))$$ and so in particular one has the implications
$$E \text{ is Ulrich }\implies E \text{  is almost Ulrich }\implies E \text{ is weakly Ulrich}.$$  
\section{Weakly Ulrich bundles on surfaces  in $\P^3$}
\newcommand{\E}{(F_*\O_X(d-3))(1)}
For a smooth, projective variety $X$ let $\omega_X$ be its canonical bundle i.e. the dualizing line bundle of $X$. \emph{In this section I  write $\deg(X)=d$}.
\bthm\label{th:main}
Let $X\subset\P^3$ be a smooth projective surface of degree $\deg(X)=d$ defined over an algebraically closed field of characteristic $p>0$.
Let $F:X\to X$ be the absolute Frobenius morphism of $X$. Let $E=\E$ and assume
$d-3 < p$.
Then $E$ is an almost Ulrich bundle on $X$.
\ethm

\blem\label{le:H1-control} 
In the above notation $H^1(X,\O_X(m))=0$  for all $m\in\Z$.
\elem
\bp 
This is standard as line bundles on $\P^3$ have no cohomology in dimensions one and two (see \cite[Chap III, Theorem 5.1]{hartshorne-algebraic})   and as $X$ is a hypersurface (see \cite[Chap III, Exercise 5.5(c)]{hartshorne-algebraic}).
\ep
\blem\label{le:omega}
With notation as above  one has $$\omega_X=\O_X(d-4).$$
\elem 
\bp 
This is standard as $X\subset \P^3$ is smooth hypersurface of degree $d$ (see \cite[Chap II, Prop. 8.20]{hartshorne-algebraic}).
\ep
\blem\label{le:H1-control2} 
For $E=\E$ one has $H^1(E(m))=0$  for all $m\in\Z$ i.e. $E$ is an ACM bundle on $X$.
\elem
\bp 
It was observed more generally in (\cite{joshi17}) that if $V$ is an ACM bundle then $F_*(V)$ is an ACM bundle. But let me provide a proof of what is asserted here for completeness. Observe that, by the projection formula (\cite[Chap II, Exercise 5.1(d)]{hartshorne-algebraic}) for $F$, one has 
\be 
E(m)=F_*(\O_X(d-3))(1+m)=F_*(\O_X(d-3+p(1+m))).
\ee As Frobenius is a finite flat morphism, one sees that, for all $i\geq 0$, one has
\be\label{eq:frobenius-cohomology} 
H^i(E(m))=H^i(F_*(\O_X(d-3+p(1+m))))=H^i(\O_X(d-3+p(1+m))).
\ee
By Lemma~\ref{le:H1-control} this is zero for $i=1$. This proves the claim.
\ep
\bp 
To prove Theorem~\ref{th:main} I have to prove the following vanishing holds:
\benum[label={{\bf(V.\arabic{*})}}] 
\item\label{v:2} $H^2(E(-m))=0$ for all $m\leq \dim(X)-1=1$,
\item\label{v:1}  $H^1(E(-m))=0$ for all $m\leq 0$ and  for $m\geq 3$,
\item\label{v:0} $H^0(E(-m))=0$ for $m\geq 2$.
\eenum
Lemma~\ref{le:H1-control2} establishes the vanishing of $H^1(E(m))$  for all $m\in\Z$ and hence also in the required range for  \ref{v:1}. Let me prove \ref{v:2}. By \eqref{eq:frobenius-cohomology} one has to prove that $H^2(\O_X(d-3+p(1-m)))=0$. If $m=1$ then one has (for all $d\geq 1$)
\be 
H^2(\O_X(d-3))=H^0(\O_X(-(d-3)\tensor \omega_X)=H^0(\O_X(-d+3+d-4))=H^0(\O_X(-1))=0
\ee
For $m=0$ one has $H^2(\O_X(d-3+p))=0$ provided, that by Serre duality, $$H^0(\O_X(-d+3-p+d-4))=H^0(\O_X(-1-p))=0$$ which is clear.

So now assume $m<0$ which says $1-m>0$. Then $H^2(\O_X(d-3+p(1-m))$ is dual to 
$$H^0(\O_X(-d+3-p(1-m)+d-4))=H^0(\O_X(-1-p(1-m)))=0 \text{ as } m<0.$$
This proves \ref{v:2}.

Let me prove \ref{v:0} which says $H^0(\O_X(d-3+p(1-m)))=0$ for $m\geq 2$. Observe that $m\geq 2$ says that $1-m\leq -1$; and, hence to prove the required vanishing, it suffices to check that
$$d-3+p(1-m)\leq d-3-p <0.$$
This is true by my hypothesis $d-3<p$. Hence \ref{v:0} is proved.
\ep
\section{Weakly Ulrich bundles on a class of surfaces of general type}
\renewcommand{\E}{(F_*\omega_X)(1)}
The method introduced above can be applied to some other situations. Here is a variant of the method sketched above.
\bthm\label{th:gen-type} 
Suppose $X$ is a smooth, projective and minimal surface of general type and that $\omega_X$ very ample and provides an embedding $X\into\P^n$ and $X$ equipped with $\O_X(1)=\omega_X$ given by this embedding. Assume
\benum 
\item $H^1(X,\O_X)=0$,
\item $p>2$.
\eenum
Then $E=F_*(\omega_X)(1)$ is an almost Ulrich bundle on $X$.
\ethm 
\bp 
It is important to note that I do not assume that Kodaira vanishing holds for $X$. I claim that $H^1(\omega^m_X)=0$ for all $m\in\Z$. If $m=0$ this is $H^1(\O_X)=0$ which is my hypothesis. If $m<0$ then one has by (\cite[Theorem~1.7]{ekedahl88}) which says, as $p>2$, that $H^1(\omega^m_X)=0$ for $m<0$ as $\omega_X$ is ample (hence nef and big). If $m>0$ then $H^1(\omega^m_X)$ is dual to $H^1(\omega_X^{1-m})$. If $m=1$ this is $H^1(\O_X)=0$ so one can assume $m\geq 2$ and hence $1-m<0$ and again one is done by (\cite[Theorem 1.7]{ekedahl88}). Thus the assertion is proved. 

I claim that $H^1(\omega_X(m))=0$ for all $m\in \Z$. This is immediate from the above vanishing as $\O_X(1)=\omega_X$. 
Next I claim that $H^1(\E(m))=0$ for all $m\in\Z$. Indeed $H^1(\E(m))=H^1(\omega_X(p(m+1)))=H^1(\omega_X^{1+p(1+m)})=0$ for all $m\in\Z$ as $H^1(\omega^m_X)=0$ for all $m\in\Z$.

Thus to prove the theorem it remains to prove $H^2(E(-m))=0$ for $m\leq 1$ and $H^0(E(-m))=0$ for $m\geq 2$.
 
Consider $H^2(E(-m))$ for $m\leq 1$. If $m=1$ then $H^2(E(-1))=H^2(F_*\omega_X)=H^2(\omega_X^p)=0$ as  its dual is $H^0(\omega^{1-p}_X)=0$ by ampleness. If $m\leq 0$ then $H^2(E(-m))=H^2(\omega_X^{1+p(1-m)})$ which is dual to $H^0(\omega_X^{-p(1-m)})=0$ as $1-m\geq 1$ for $m\leq 0$ and by ampleness. 

So it remains to prove that $H^0(E(-m))=0$ for $m\geq 1$. If $m\geq 2$ this is $H^0(\omega_X^{1+p(1-m)})$ with $1+p(1-m)\leq -1$ and so one is done by ampleness. So $E$ is weakly Ulrich and ACM hence $E$ is an almost Ulrich bundle as claimed.
\ep

\brem 
Let me remark that surfaces which satisfy the hypothesis of Theorem~\ref{th:gen-type} exist. For example a quintic surface $X\subset\P^3$ satisfies all the hypothesis of the theorem.
\erem
\section{Weakly Ulrich bundles on Frobenius split varieties of dimension at most three}
Let $X/k$ be a smooth, projective  variety of dimension $\dim(X)=q$ over an algebraically closed field $k$ of characteristic $p>0$. Let $F:X\to X$ be the absolute Frobenius morphism of $X$. Then one has an exact sequence (which can be taken to be the defining the vector bundle $B^1_X$):
\be\label{eq:b1-def} 
0\to \O_X \to F_*(\O_X) \to B^1_X \to 0,
\ee
and $B^1_X$ is locally free of rank $p^q-1$. Following (\cite{mehta85}) one says that $X$ is Frobenius split if \eqref{eq:b1-def} splits as a sequence of $\O_X$-modules. One has the following (for a proof see (\cite{joshi03})):
\be\label{eq:b1-cohom} 
H^i(X,B^1_X)=0 \text{ for all }i\geq 0.
\ee

Let us note that if $X$ is a smooth, projective and Frobenius split curve then $X\isom \P^1$ or $X$ is an ordinary elliptic curve. In this case existence of Ulrich bundles is elementary and can be treated as in (\cite{joshi17}) where I was primarily concerned with the case of curves of genus at least two.

For Frobenius split varieties one has  
the following vanishing theorem due to (\cite{mehta85}). This result is stronger than Kodaira vanishing theorem (note that this form of the result is not stated in loc. cit.).
\bthm\label{th:mehta}[\cite{mehta85}]
Let $X/k$ be a smooth, projective and Frobenius split variety and suppose $L$ is an ample line bundle on $X$. Then
\benum[label={{\bf(\arabic{*})}}]
\item  $H^i(X,L)=0$ for all $i>0$, and
\item $H^i(X,L^{-1})=0$ for all $i<\dim(X)$.
\eenum
\ethm 

\bthm\label{th:calabi-yau}
Let $X$ be a Frobenius split variety of dimension $\dim(X)=q\geq 2$ with an embedding $X\into \P^n$ and equipped with $\O_X(1)$ given by this embedding.  Let $E=B^1_X(q-1)$. Then
\benum[label={{\bf(\arabic{*})}}]
\item\label{th:calabi-yau1} for all $q\geq 2$, $E$ is an almost Ulrich bundle on $X$ if and only if $H^0(B^1_X(q-1-m))=0$ for $2\leq m\leq q-2$.
\item\label{th:calabi-yau2}  If $q=2,3$ then $E$ is an almost Ulrich bundle on $X$.
\item\label{th:calabi-yau3}  If $q=1$ then $B^1_X(1)$ is an Ulrich bundle on $X$.
\eenum
\ethm
\bp 
In (\cite{joshi17}) it was shown that under the above hypothesis of \ref{th:calabi-yau1} that $B^1_X$ is an ACM bundle on $X$ and hence so is $E=B^1_X(q-1)$. Thus it follows that $E=B^1_X(q-1)$ has no cohomology in dimensions $1\leq i\leq q-1=\dim(X)-1$. 
Now let me prove Theorem~\ref{th:calabi-yau}\ref{th:calabi-yau1}. To prove that $E$ is weakly Ulrich bundle one has to deal with cohomologies in dimensions $i=0,q$. So we have to prove:
\benum[label={{\bf(C.\arabic{*})}}]
\item\label{c:d} $H^q(E(-m))=0$ for all $m\leq q-1$,
\item\label{c:0} $H^0(E(-m))=0$ for all $m\geq 2$.
\eenum
Let me begin by proving \ref{c:d}. As \eqref{eq:b1-def} splits, cohomology of $B^1_X$ and all its twists are direct summands of the cohomology of $F_*(\O_X)$ and its twists. By the projection formula (\cite[Chap II, Exercise 5.1(d)]{hartshorne-algebraic}) $$F_*(\O_X)(q-1-m)=F_*(\O_X(p(q-1-m))).$$ Hence, to prove that $H^q(E(-m))=0$ in the asserted range, it suffices to prove the vanishing of $$H^q(F_*(\O_X(p(q-1-m))))=H^q(\O_X(p(q-1-m))).$$ First let us deal with the exceptional case  $q-1-m=0$ i.e.  $m=q-1$. In this case, by \eqref{eq:b1-cohom}, one has $H^q(E(-(q-1)))=H^q(B^1_X)=0$. So one can assume that $m<q-1$. In this case $\O_X(p(q-1-m))$ is ample (for all $m<q-1$); and, as $X$ is Frobenius split, so by Theorem~\ref{th:mehta},  one has  that $$H^q(\O_X(p(q-1-m)))=0 \qquad \text{ for all } m<q-1.$$ 

So now it remains to prove \ref{c:0}. As $H^0(E(-m))$ is a direct summand of $H^0(\O_X(p(q-1-m)))$, one sees that if $m\geq q$ then $q-1-m\leq q-1-q<0$; and, so one can appeal to Kodaira vanishing theorem of (\cite{mehta85}). So it remains to deal with the range $2\leq m\leq q-1$. If $m=q-1$ then $q-1-m=0$ and so one has to prove that $$H^0(E(-(q-1)))=H^0(B^1_X)=0$$ which is the case by \eqref{eq:b1-cohom}. For $2\leq m\leq q-2$ the required vanishing is our hypothesis. This completes the proof of Theorem~\ref{th:calabi-yau}\ref{th:calabi-yau1}. 

Observe that Theorem~\ref{th:calabi-yau}\ref{th:calabi-yau1} $\implies$ Theorem~\ref{th:calabi-yau}\ref{th:calabi-yau2} since for $q=2,3$ the condition $2\leq m\leq q-2$ is vacuously true. Hence one sees that for $q=2,3$  the bundle $E=B^1_X(q-1)$ is a weakly Ulrich bundle; and, for all $q\geq 2$, $E$ is also an ACM bundle on $X$. 

Finally note that if $q=1$ then the fact that $B^1_X(1)$ is an Ulrich bundle on $X$ is immediate from (\cite{joshi17}). Note that  any smooth, projective curve $X$ is Frobenius split curve if and only if $X=\P^1$ (which is also ordinary) or $X$ is an ordinary elliptic curve. For a comparison of ordinary and Frobenius split varieties reader should consult (\cite{joshi03}).
\ep

This result  includes the following special cases:
\bcor\label{cor:calabi-yau}
Let $X/k$ be a smooth, projective Frobenius split Calabi-Yau variety with $2\leq\dim(X)\leq 3$ (in particular any ordinary K3 surface)  then $E=B^1_X(q-1)$ is an almost Ulrich bundle on $X$. 
\ecor

\bcor\label{cor:abelian} 
Let $X/k$ be a Frobenius split (i.e. ordinary) abelian surface or an abelian threefold then $E=B^1_X(q-1)$ is an almost Ulrich bundle on $X$.
\ecor

\bcor 
Let $X/k$ be a smooth, projective, Frobenius split and Fano variety dimension with $2\leq\dim(X)\leq 3$. 
Then $E=B^1_X(q-1)$ is an almost Ulrich bundle on $X$.
\ecor
Note that $\P^n$ is a Frobenius split Fano variety in any dimension $n\geq 1$. Hence the above corollary also applies (quite trivially) to any smooth embeddings  $X=\P^2\into\P^n$ and $X=\P^3\into \P^n$.

In particular as any Fano hypersurface in $\P^n$ is Frobenius split by (\cite{fedder87}) one has the following:
\bcor\label{cor:fano} 
Let $X/k$ be any smooth, Fano hypersurface of dimension with $2\leq\dim(X)\leq 3$ in $\P^n$  (so $n=3,4$ and degree of $X$ is at most $n$) then $E=B^1_X(q-1)$ is an almost Ulrich bundle on $X$. 
\ecor

\brem\label{re:almost-ulrich+epsilon} 
Let me remark that in Theorem~\ref{th:calabi-yau} the bundle $E=B^1_X(q-1)$ also satisfies $H^0(E(-1))=0$ so $E$ is an almost Ulrich bundle which fails to be Ulrich only because $H^q(E(-q))$ may not be zero. 
\erem

Note that the bundle $E=B_X^1(q-1)$ is intrinsic to $X$: its dependence on the embedding $X\into\P^n$ comes only from the presence of the twist by $\O_X(q-1)$. So these corollaries together with (\cite[Section 2, Section 3]{eisenbud03}) imply the following (see loc. cit. for the definition of the Chow form of a variety).
\bthm\label{th:chow-form}
Let $k$ be an algebraically closed field of characteristic $p>0$. Let $X/k$  be a smooth, projective, Frobenius split variety (in particular a Frobenius split Calabi-Yau, Abelian or Fano variety) of dimension at most three equipped with a projective embedding $X\into\P^n$. Then the Chow form of $X$ with respect to this embedding is the support of a single intrinsic determinantal equation.
\ethm

\bibliographystyle{plainnat}
\bibliography{../ulrich/ulrich.bib,../../master/joshi.bib,../../master/master6.bib}
\end{document}